\documentclass[12pt,reqno]{amsart}
\usepackage{amsmath, amssymb, latexsym, amscd, amsthm,amsfonts,amstext}
\usepackage[mathscr]{eucal}
\usepackage{xspace}
\usepackage{array, tabularx, longtable}
\usepackage{multicol,color}
\usepackage{indentfirst}
\usepackage{graphics}
\usepackage[colorlinks=true]{hyperref}
\theoremstyle{plain}
\newtheorem{theorem}{Theorem}

\newtheorem{lemma}[theorem]{Lemma}
\newtheorem{proposition}[theorem]{Proposition}
\theoremstyle{remark}

\newtheorem{claim}{Claim}
\newtheorem{remark}{Remark}
\theoremstyle{definition}
\newtheorem{define}{Definition}

\begin{document}
\date{} 

\title [Modified defect relations of the Gauss map on annular ends ]{Modified defect relations of the Gauss map of complete minimal surfaces on annular ends}

\author{Pham Hoang Ha and Nguyen Hoang Trang}


\begin{abstract} 
In this article, we study the modified defect relations of the Gauss map of complete minimal surfaces in $\mathbb R^3$ and $ \mathbb R^4$ on annular ends. We obtain results which are similar to the ones obtained by Fujimoto~[J. Differential Geometry \textbf{29} (1989), 245-262] for (the whole) complete minimal surfaces. We thus give some improvements of the previous results for the Gauss maps of complete minimal surfaces restricted on annular ends.
\end{abstract}
\maketitle
\tableofcontents
\section{Introduction}
Let $M$ be a non-flat minimal surface in $\mathbb R^3,$ or more precisely, a non-flat connected oriented minimal surface in $\mathbb R^3.$ By definition, the Gauss map $G$ of $M$ is a map which maps each point $p \in M$ to the unit normal vector $G(p) \in S^2$ of $M$ at $p.$ Instead of $G,$ we study the map $g:= \pi\circ G: M \rightarrow  \mathbb P^1(\mathbb C)$, where  $\pi : S^2 \rightarrow \mathbb P^1(\mathbb C)$
is the stereo\-graphic projection.  By associating a holomorphic local coordinate $z=u+\sqrt{-1}v$ with each positive isothermal coordinate system $(u, v),$ $M$ is considered as an open Riemann surface with a conformal metric $ds^2$ and by the assumption of minimality of $M,$ $g$ is a meromorphic function on $M.$ 

In 1988, H. Fujimoto (\cite{Fu1}) proved Nirenberg's conjecture that if $M$ is a complete non-flat minimal surface in $\mathbb R^3,$ then its Gauss map can omit at most 4 points, and the bound is sharp. After that, he (\cite{Fu2}) also extended that result by introduced the modified defect. Beside that, Ru (\cite{Ru2}), Ros (\cite{Ro}), Osserman-Ru (\cite{OR}) and the first author (\cite{H}) also improved Fujimoto's result by using the ramifications of the Gauss map or estimating the Gauss curvature of minimal surfaces.

On the other hand, in 1991, S. J. Kao (\cite{Kao}) used the ideas of Fujimoto~(\cite{Fu1}) to show that the Gauss map of an end of a non-flat complete minimal surface in $\mathbb R^3$ that is conformally an annulus $\{ z  : 0 < 1/r < |z| < r \}$ must also assume every value, with at most 4 exceptions. In 2007, L. Jin and M. Ru (\cite{JR}) extended Kao's result to minimal surfaces in $\mathbb R^m.$ Recently, Dethloff-Ha (\cite{DH}), Dethloff-Ha-Thoan (\cite{DHT}) gave some improvements for the results of Kao and Jin-Ru by estimating the ramifications of the Gauss map.

A natural question is whether a result for the modified defect relations of the Gauss map still holds on an annular end of  a non-flat complete minimal surface in $\mathbb R^3$. In this paper we give an affirmative answer for this question. Moreover, we also prove some analogue results for the Gauss map of minimal surface in $\mathbb R^4.$ We thus give some improvements of the results of Kao (\cite{Kao}) and Dethloff-Ha (\cite{DH}).
\section{Statements of the main results}

 Firstly, we recall some definitions and remarks on the modified defect which was introduced by Fujimoto (\cite{Fu2}).

Let $M$ be an open Riemann surface and $f$ a nonconstant holomorphic map of $M$ into $\mathbb P^1(\mathbb C).$ Assume that $f$ has reduced representation $f=(f_0: f_1).$ Set $||f|| = (|f_0|^2 + |f_1|^2)^{1/2}$ and, for each $\alpha = (a_0: a_1) \in \mathbb P^1(\mathbb C)$ with $|a_0|^2 + |a_1|^2 = 1,$ we define $F_\alpha := a_1f_0 - a_0f_1.$
\begin{define} \label{d1}
We define the $S-$defect of $\alpha$ for $f$ by
$$\delta^S_f(\alpha):= 1 - \inf \{\eta \geq 0; \eta \text{ satisfies condition $(*)_S$}\}.$$
Here, condition $(*)_S$ means that there exists a $ [-\infty, \infty)-$valued continuous subhamornic function $u\  (\not\equiv - \infty)$ on $M$ satisfying the following conditions:\\
(C1) $e^{u} \leq ||f||^{\eta},$ \\
(C2) for each $\xi \in f^{-1}(\alpha) ,$ there exists the limit\\
$$\lim_{z \rightarrow \xi}(u(z) - \log|z-\xi|) \in [-\infty, \infty) ,$$
where $z$ is holomorphic local coordinate around $\xi.$
\end{define}
\begin{remark}
We always have that $\eta = 1$ satifies condition $(*)_S$ when we consider $u = |F_\alpha|.$
\end{remark}
\begin{define}
\label{d2}
We define the $H-$defect of $\alpha$ for $f$ by
$$\delta^H_f(\alpha):= 1 - \inf \{\eta \geq 0; \eta \text{ satisfies condition $(*)_H$}\}.$$
Here, condition $(*)_H$ means that there exists a $ [-\infty, \infty)-$valued continuous subhamornic function $u$ on $M$ which is hamornic on $M - f^{-1}(\alpha)$ and satisfies the  conditions (C1) and (C2).
\end{define}
\begin{define} \label{d3}
We define the $O-$defect of $\alpha$ for $f$ by
$$\delta^O_f(\alpha):= 1 - \inf \{\ \dfrac{1}{m}; \ \text{ $F_\alpha$ has no zero of order less than $m$.}\}$$
\end{define}
\begin{remark}
We always have $0\leq \delta^O_f(\alpha ) \leq \delta^H_f(\alpha)\leq \delta^S_f(\alpha)\leq 1.$
\end{remark}
On the other hand we have
\begin{define}
$f$ is called to be {\it  ramified over a point} $\alpha \in  \mathbb P^{1}(\mathbb C)$ {\it with multiplicity at least} $m$ if all the zeros of the function $ F_\alpha$ have orders at least $m.$ If the image of $f$ omits $\alpha,$ we will say that $f$ is {\it ramified over a with multiplicity }$\infty.$
\end{define}
\begin{remark} \label{rm2}
If  $f$ ramified over a point $\alpha \in  \mathbb P^{1}(\mathbb C)$ with multiplicity at least $m,$ then $\delta^H_f(\alpha) \geq \delta^O_f(\alpha ) \geq 1 - \dfrac{1}{m}.$ In particular, if $f^{-1}(\alpha) = \emptyset,$ then $\delta^O_f(\alpha ) = 1.$
\end{remark}
The main purpose of this article is to prove some results on the modified defect relations of the Gauss map of minimal surfaces on anular ends. We first prove the following theorem.
\begin{theorem} \label{T1}
Let $M$ be a non-flat complete minimal surface in $\mathbb R^3$ with the Gauss map $g$ and let $A \subset M$ be an
annular end of $M$ which is conformal to $\{z :  0 < 1/r < |z| < r\}$, where $z$
is a conformal coordinate. For arbitrary $q$ distinct points $a^1, ..., a^q$ in $\mathbb P^1(\mathbb C),$ then
\begin{equation*}
 \sum_{j=1}^q\delta^H_g(a^j) \leq 4.
\end{equation*}  
\end{theorem}
By Remark \ref{rm2} and Theorem \ref{T1} gives the following generalization of the following theorem of Dethloff-Ha (\cite{DH}).\\
{\bf Theorem A.} \ {\it Let $M$ be a non-flat complete minimal surface in $\mathbb R^3$ and let $A \subset M$ be an annular end of $M$ which is conformal to $\{z  :  0 < 1/r < |z| < r\}$, where $z$ is a conformal coordinate. If there are $q\  (q > 4)$ distinct points $a^1, ... , a^q \in \mathbb P^1(\mathbb C)$ such that the restriction of the Gauss map $g$ of $M$ to $A$ is ramified over $a^j$ with multiplicity at least $m_j$ for each $j$,   then 
\begin{equation} \label{1} \sum_{j=1}^q(1 - \frac{1}{m_j})\leq 4.
\end{equation}  
Moreover, (\ref{1})  still holds if we replace, for all $j=1,...,q$,  $m_j$ by the limit inferior
of the orders of the zeros of the function $ G_{a^j}$ on $A$ and in particular 
by $\infty$ if $g$ takes the value $a^j$ only a finite number of times on $A$.}

Moreover, we also would like to consider  the Gauss map of complete minimal surfaces $M$ immersed in $\mathbb R^4,$ this case has been investigated by various authors (see, for example Osserman (\cite{Os}), Chen (\cite{Chen}), Fujimoto (\cite{Fu2}) and Dethloff-Ha (\cite{DH})). In this case, the Gauss map $g$ of $M$ may be identified with a pair of meromorphic functions $g=(g^1, g^2).$ For the last purpose of this article, we shall prove the following result of modified defect relations of the Gauss map restricted on an annular end. 
\begin{theorem}\label{T2}
Suppose that $M$ is a complete non-flat minimal surface in $\mathbb R^4$ and $g=(g^1, g^2)$ is the Gauss map of $M.$ Let A be an annular end of $M$ which is conformal to $\{z  :  0 < 1/r < |z| < r\}$, where $z$
is a conformal coordinate. Let $a^{11}, ... , a^{1q_1}, a^{21}, ... , a^{2q_2}$ be $ q_1 + q_2\ (q_1, q_2 > 2)$
distinct points in $ \mathbb P^1(\mathbb C).$ \\
\indent (i) In the case $g^l\not\equiv constant\  (l=1,2),$  then
$ \sum_{j=1}^{q_1}\delta^H_{g^1}(a^{1j})\leq 2,$\\ or $ \sum_{j=1}^{q_2}\delta^H_{g^2}(a^{2j})\leq 2,$ or
$$\dfrac{1}{\sum_{j=1}^{q_1}\delta^H_{g^1}(a^{1j}) - 2} + \dfrac{1}{\sum_{j=1}^{q_2}\delta^H_{g^2}(a^{2j}) - 2} \geq 1.$$
\indent (ii) In the case where one of $g^1$ and $g^2$ is constant, say $g^2\equiv constant, $ we have the following
$$ \sum_{j=1}^{q_1}\delta^H_{g^1}(a^{1j}) \leq 3.$$
\end{theorem}
Combining with Remark \ref{rm2}, we get directly the following theorem of Dethloff-Ha (\cite{DH}) from Theorem \ref{T2}.\\
{\bf Theorem B.}\ {\it Suppose that $M$ is a complete non-flat minimal surface in $\mathbb R^4$ and $g=(g^1, g^2)$ is the Gauss map of $M.$ Let A be an annular end of $M$ which is conformal to $\{z  :  0 < 1/r < |z| < r\}$, where $z$
is a conformal coordinate. Let $a^{11}, ... , a^{1q_1}$ respectively $a^{21}, ... , a^{2q_2}$ be $ q_1\ (q_1>2)$ respectively  $q_2\ (q_2>2)$ 
distinct points in $ \mathbb P^1(\mathbb C).$ \\
\indent (i) In the case $g^l\not\equiv constant\  (l=1,2),$ if $g^l$ is ramified over $a^{lj}$ with multiplicity at least $m_{lj}$ for each $j\ (l=1,2)$ on $A,$  then\\
$ \gamma_1=\sum_{j=1}^{q_1}(1 - \frac{1}{m_{1j}})\leq 2,$ or $ \gamma_2=\sum_{j=1}^{q_2}(1 - \frac{1}{m_{2j}})\leq 2,$ or\\
$$\dfrac{1}{\gamma_1 - 2} + \dfrac{1}{\gamma_2 - 2} \geq 1.$$
\indent (ii) In the case where  $g^1$ or $g^2$ is constant, say $g^2\equiv constant, $ if $g^1$  is ramified over $a^{1j}$ with multiplicity at least $m_{1j}$ for each $j$ on $A,$ we have the following :
$$ \gamma_1=\sum_{j=1}^{q_1}(1 - \frac{1}{m_{1j}})\leq 3.$$
Moreover, the results  still hold if we replace, for all $a^{lj}$ $(j=1,...,q_l; \:l=1,2)$  the $m_{lj}$ by the limit inferior
of the orders of the zeros of the function $G_{a^{lj}}$ on $A$ and in particular 
by $\infty$ if $g^l$ takes the value $a^{lj}$ only a finite number of times on $A$.}

The main idea to prove our theorems is to construct a pseudo-metric with negative curvature with conditions of modified defect relations on an annular end, which is a refinement of  the ideas in Fujimoto (\cite{Fu2}) and Dethloff-Ha (\cite{DH}). After that we use  arguments similar to those used by Kao (\cite{Kao}) and by Dethloff-Ha (\cite{DH}) to finish the proofs. 

\section{Auxiliary  lemmas}
Let $f$ be a nonconstant holomorphic map of a disk $\Delta_R := \{z \in \mathbb C: |z| < R \}$ into $\mathbb P^1(\mathbb C),$ where $0 < R < \infty.$ Take a reduced representation $f = (f_0: f_1)$ on $\Delta_R$ and define
$$||f||:= (|f_0|^2 + |f_1|^2)^{1/2}, W(f_0,f_1) =W_z(f_0, f_1):= f_0f_1' - f_1f_0',$$
where the derivatives are taken with respect to the variable $z$.
Let $a^j\: ( 1\leq j \leq q )$ be $q$ distinct points in $\mathbb P^1(\mathbb C).$ We may assume $a^j=(a^j_0: a^j_1 )$ with $|a^j_0|^2 + |a^j_1|^2 = 1$ $( 1\leq j \leq q ),$ and set
$$F_j:=a^j_0f_1 - a^j_1f_0 \ (1\leq j \leq q).$$
\begin{proposition}(\cite[Proposition 2.1.6 and 2.1.7]{Fu3})\label{W}.\\
a) If $\xi$ is another local coordinate, then 
$W_{\xi}(f_0,f_1)= W_z(f_0,f_1) \cdot (\frac{dz}{d\xi})$.\\
b)  $W(f_0,f_1) \not\equiv 0$ (iff $f$is nonconstant).
\end{proposition}
We now consider $[-\infty, \infty)-$valued continuous subhamornic functions $u_j (\not\equiv - \infty)$ on $\Delta_R$ and nonnegative numbers $\eta_j (1 \leq j \leq q)$ satisfying the conditions:\\
(D0) $\gamma := q - 2 - (\sum_{j=1}^q \eta_j) > 0$\\
(D1) $e^{u_j} \leq ||f||^{\eta_j}$ for $j=1, ..., q,$\\
(D2) for each $\xi \in f^{-1}(\alpha_j) \ ( 1\leq j \leq q ),$ there exists the limit
$$\lim_{z \rightarrow \xi}(u_j(z) - \log|z-\xi|) \in [-\infty, \infty) .$$
\begin{lemma}(\cite{Fu2})
\label{l1}
Let $u_1, \cdots , u_q$ be continuous subharmonic functions on $\Delta_R$, and $\eta_1, \cdots , \eta_q$ nonnegative constants which satisfy the conditions (D0)-(D2). Then, for each $\sigma$ with $0< q\delta < \gamma,$ there exists a positive constants $C$ such that 
$$ \dfrac{||f||^{\gamma-q\delta}e^{\sum_{j=1}^qu_j}|W(f_0,f_1)|}{\Pi_{j=1}^q|F_j|^{1-\delta}} \leq C\dfrac{2R}{R^2-|z|^2}.$$
\end{lemma}

\begin{lemma}(\cite[Lemma 1.6.7]{Fu3}). \label{L5}
Let $d\sigma^2$ be a conformal flat metric on an open Riemann surface $M$. Then for every point $p \in M$, there is a holomorphic and locally biholomorphic map $\Phi$ of
a disk (possibly with radius $\infty$)  $\Delta_{R_0} := \{w : |w|<R_0 \}$ $(0<R_0 \leq \infty )$ onto an open neighborhood of $p$ with $\Phi (0) = p$ such that $\Phi$ is a local isometry, namely the pull-back 
$\Phi^*(d\sigma^2)$ is equal to the standard (flat) metric on $\Delta_{R_0}$, and for some point $a_0$ with $|a_0|=1$, the $\Phi$-image of the curve 
$$L_{a_0} : w:= a_0 \cdot s \; (0 \leq s < R_0)$$
is divergent in $M$ (i.e. for any compact set $K \subset M$, there exists an $s_0<R_0$
such that the $\Phi$-image of the curve $L_{a_0} : w:= a_0 \cdot s \; (s_0 \leq s < R_0)$
does not intersect $K$).
\end{lemma}
\section{The proof of  Theorem \ref{T1}}
\begin{proof}
\indent For convenience of the reader, we first recall some notations on the Gauss map of minimal surfaces in $\mathbb R^3.$
\indent Let $x =( x_1, x_2, x_3) : M \rightarrow \mathbb R^3$ be a non-flat complete minimal surface and $g: M \rightarrow  \mathbb P^1(\mathbb C)$ its Gauss map.  Let $z$ be a local holomorphic coordinate. Set $\phi_i := \partial x_i / \partial z \ (i = 1, 2, 3 )$ and $\phi:= \phi_1-\sqrt{-1}\phi_2.$ Then, the (classical) Gauss map $g: M \rightarrow \mathbb P^1(\mathbb C)$ is given by  $$g=\dfrac{\phi_3}{\phi_1 - \sqrt{-1}\phi_2},$$
and the metric on $M$ induced from $\mathbb R^3$ is given by
\begin{equation}\label{eq:41}
ds^2= |\phi|^2(1 + |g|^2)^2|dz|^2  \text{ (see Fujimoto (\cite{Fu3}))}.
\end{equation}
We remark that although the $\phi_i$, $(i=1,2,3)$ and $\phi$ depend on $z$, $g$ and $ds^2$ do not.
Next we take a reduced representation $g = (g_0 : g_1)$ on $M$ and set $||g|| = (|g_0|^2 +|g_1|^2)^{1/2}.$ Then we can rewrite 
\begin{equation}\label{eq:42}
ds^2 = |h|^2||g||^4|dz|^2\,, 
\end{equation}
where $h:= \phi/g_0^2.$ In particular, $h$ is a holomorphic map without zeros. We remark that  $h$ depends on $z$, however, the reduced representation $g=(g_0:g_1)$ is globally defined on $M$ and independent of $z$.
Finally we observe that by the assumption that $M$ is not flat,  $g$ is not constant.

Now the proof of Theorem \ref{T1} will be given in three steps :

{\bf Step 1:} Let $a^j\:( 1\leq j \leq q )$ be $q$ distinct points in $\mathbb P^1(\mathbb C).$ We may assume $a^j=(a^j_0: a^j_1 )$ with $|a^j_0|^2 + |a^j_1|^2 = 1$ $( 1\leq j \leq q ),$ and we set $G_j:=a^j_0g_1 - a^j_1g_0 \ (1\leq j \leq q)$ for the reduced representation $g = (g_0 : g_1)$ of the Gauss map. Let $A \subset M$ be an annular end of $M,$ that is, $A = \{z  :  0 < 1/r < |z| < r < \infty  \},$ where $z$ is a (global) conformal coordinate of $A$. All the $\delta^H_g(a^j)$ are increasing if we only consider the Gauss map $g$ which takes on a subset $B \subset A$. So without loss of generality we may prove our theorem only
 on a sub-annular end, i.e. a subset $ A_t :=\{z  :  0 < t \leq |z| < r < \infty  \} \subset A$ 
 with some $t$ such that $1/r < t<r$. (We trivially observe that for $c:=tr>1$, $s:= r/\sqrt{c}$, $\xi := z/ \sqrt{c}$, we have $A_t = \{\xi  :  0 < 1/s \leq |\xi| < s < \infty  \}$). 

 Now on the annular end $A = \{z  :  0 < 1/r \leq |z| < r < \infty  \},$ we may assume that 
 \begin{equation} \label{ass2} \sum_{j=1}^q\delta^H_g(a^j) > 4,
\end{equation}  
since otherwise Theorem \ref{T1} is already proved. 

By definition, there exist constants $\eta_j \geq 0 (1\leq j \leq q)$ such that $\gamma := q-2-\sum_{j=1}^q\eta_j > 2$ and continous functions $u_j (1\leq j \leq q)$ on $M$ such that each $u_j$ is harmonic on $M\setminus g^{-1}(a^j)$ and satisfies conditions (D1) and (D2). Take $\delta$ with
\begin{equation}\label{eq:44}
\dfrac{\gamma -2}{q} > \delta > \dfrac{\gamma - 2}{q +2},
\end{equation} 
and set $p = 2/ (\gamma -q\delta).$ Then  
\begin{equation} \label{3.4.1}0 < p < 1 , \ \frac{p}{1-p} > \frac{\delta p}{1-p} > 1 \ . 
\end{equation}
 Consider the  subset
$$ A_1 = A \setminus \{ z : \Pi_{j=1}^qG_j(z)\cdot W_z(g_0, g_1)(z) = 0 \}$$
of $A$. We define a new metric
\begin{equation}\label{eq:46}
d\tau^2 = |h|^{\frac{2}{1-p}}\bigg(\dfrac{\Pi_{j=1}^q|G_j|^{1-\delta}}{e^{\sum_{j=1}^qu_j}|W(g_0,g_1)|}\bigg)^{\frac{2p}{1-p}}|dz|^2
\end{equation}
on $A_1$ (where again $G_j := a^j_0g_1 - a^j_1g_0$ and $h$ is defined with respect
to the coordinate $z$ on $A_1 \subset A$ and $W(g_0,g_1) = W_z(g_0,g_1)$). 

First we observe that $d\tau$ is continuous and nowhere vanishing on $A_1.$ Indeed, $h$ is without zeros on $A$ and for each $z_0 \in A_1$ with $G_j(z_0)\cdot W(g_0,g_1)(z_0) \not= 0$ for all $j=1,...,q$,  $d\tau$ is continuous at $z_0.$ Thus, $d\tau$ is continuous and nowhere vanishing on $A_1.$ \\
Next, it is easy to see that $d\tau$ is flat.\\
By Proposition \ref{W} a) and the dependence of $h$ on $z$ and the independence of the $G_j$ of $z$, we also easily see that $d\tau$ is independent of the choice of the coordinate $z$.

The key point is to prove the following claim :
\begin{claim} \label{Cl1}
 $d\tau$ is complete on the set  $\{ z : |z| = r\}\cup \{z : W(g_0, g_1)(z) = 0\},$ i.e., the set $\{ z : |z| = r\}\cup \{z : W(g_0, g_1)(z) = 0 \}$ is at infinite distance from any interior point in $A_1.$ 
\end{claim}
 If $W(g_0, g_1)(z_0) = 0,$ then we have two cases :\\
{\it Case 1.} $G_j(z_0) = 0$ for some $j \in \{1, 2, ..., q \}.$\\
Then we have $G_i(z_0) \not= 0$ for all $i\not= j$ and $\nu_{G_j}(z_0) \geq 1.$ Changing the indices if necessary, we may assume that $g_0(z_0) \not= 0$, so also $a^j_0 \not= 0.$ So, we get 
\begin{equation} \label{order}
\nu_{W(g_0,g_1)}(z_0) = \nu_{ \dfrac{( a^j_0\frac{g_1}{g_0}-a^j_1 )'}{ a^j_0}}(z_0) =\nu_{ \dfrac{(G_j/g_0)'}{ a^j_0}}(z_0) = \nu_{G_j}(z_0) - 1.
\end{equation}
and then the function
\begin{equation}\label{eq:48}
\dfrac{e^{u_j}|W(g_0,g_1)|}{|G_j|}=(|z-z_0||\chi_j|)e^{u_j-\log|z-z_0|},
\end{equation}
where $\chi_j=\dfrac{W(g_0,g_1)}{G_j}, $ is bounded in a neighborhood of $z_0.$ \\
Thus we have 
\begin{align*}
\nu_{d\tau}(z_0) &\leq \dfrac{p}{1-p}(-\nu_{\bigg(\dfrac{e^{u_j}W(g_0,g_1)}{G_j}\bigg)}(z_0)-\delta\nu_{G_j}(z_0))\\  
& \leq \dfrac{p}{1-p}(-\delta\nu_{G_j}(z_0))
\leq  \dfrac{-p\delta}{1-p}.\\
\end{align*}
{\it Case 2.} $G_j(z_0) \not= 0$ for all $ 1\leq j \leq q.$ \\
It is easily to see that  $\nu_{d\tau}(z_0) \leq -\dfrac{p}{1-p}.$\\
So, since $0<\delta<1$, we can find a positive constant $C$ such that 
$$|d\tau| \geq  \dfrac{C}{ |z-z_0|^{\delta p/(1-p)}}|dz| $$ 
in a neighborhood of $z_0$. Combining with (\ref{3.4.1}) we thus have that $d\tau$ is complete on $\{z : W(g_0,g_1)(z) = 0\}.$ \\
\indent Now assume that $d\tau$ is not complete on $\{z: |z| = r \}.$ Then there exists $\rho: [0, 1) \rightarrow A_1,$ where $\rho (1) \in \{z : |z| = r \},$ so that $|\rho| < \infty.$ Furthermore, we may also assume that $dist(\rho (0); \{z : |z| = 1/r\}) > 2|\rho|.$ Consider a small disk $\Delta$ with center at $\rho (0).$ Since $d\tau$ is flat, $\Delta$ is isometric to an ordinary disk in the plane (cf. e.g.  Lemma \ref{L5}). Let $\Phi: \{|w| < \eta \}\rightarrow \Delta$ be this isometry. Extend $\Phi,$ as a local isometry into $A_1,$ to the largest disk $\{|w| < R\} = \Delta_R$ possible. Then $R \leq |\rho|.$ The reason that $\Phi$ cannot be extended to a larger disk is that the image goes to the outside boundary $\{z : |z| = r\}$ of $A_1$ (it cannot go to points of $A$ with $W(g_0,g_1)=0$ since we have shown already the completeness of $A_1$ with respect 
to these points). More precisely, there exists a point $w_0$ with $|w_0| =R$ so that $\Phi(\overline{0,w_0}) = \Gamma_0$ is a divergent curve on $A.$\\
The map $\Phi(w)$ is locally biholomorphic, and the metric on $\Delta_R$ induced from $ds^2$ through $\Phi$ is given by
\begin{equation} \label{3.3} \Phi^*ds^2 =  |h \circ \Phi|^2||g \circ \Phi||^4|\frac{dz}{dw}|^2|dw|^2 \ .
\end{equation}
On the other hand, $\Phi$ is isometric, so we have
$$ |dw| = |d\tau|= \bigg(\dfrac{|h| \Pi_{j=1}^q|G_j|^{(1-\delta)p}}{e^{p\sum_{j=1}^qu_j}|W(g_0,g_1)|^p}\bigg)^{\frac{1}{1-p}}|dz|$$
$$\Rightarrow |\dfrac{dw}{dz}|^{1-p} = \dfrac{|h| \Pi_{j=1}^q|G_j|^{(1-\delta)p}}{e^{p\sum_{j=1}^qu_j}|W(g_0,g_1)|^p}.$$
Set $f:= g(\Phi), f_0 := g_0(\Phi), f_1 := g_1(\Phi), u_j := u_j(\Phi)$ and $F_j:= G_j(\Phi).$  Since
$$W_w(f_0, f_1) = (W_z(g_0, g_1) \circ \Phi)\frac{dz}{dw}, $$
we obtain
\begin{equation} \label{3.4} |\dfrac{dz}{dw}| = \dfrac{e^{p\sum_{j=1}^qu_j}|W(f_0,f_1)|^p}{|h(\Phi)| \Pi_{j=1}^q|F_j|^{(1-\delta)p}}.
\end{equation}
By (\ref{3.3}) and (\ref{3.4}) and by definition of $p$, therefore, we get
\begin{align*}
\Phi^*ds^2& = \bigg( \dfrac{||f||^2e^{p\sum_{j=1}^qu_j}|W(f_0,f_1)|^p}{ \Pi_{j=1}^q|F_j|^{(1-\delta)p}}\bigg)^2|dw|^2\\
&= \bigg( \dfrac{||f||^{\gamma-q\delta}e^{\sum_{j=1}^qu_j}|W(f_0,f_1)|}{ \Pi_{j=1}^q|F_j|^{1-\delta}}\bigg)^{2p}|dw|^2.
\end{align*}
Using the Lemma \ref{l1}, we obtain
$$\Phi^*ds^2 \leqslant C^{2p}.(\dfrac{2R}{R^2 -|w|^2})^{2p}|dw|^2.$$
 Since $0 < p < 1,$ it then follows that 
$$d_{\Gamma_0} \leqslant \int_{\Gamma_0}ds = \int_{\overline{0,w_0}}\Phi^*ds \leqslant C^p. \int_0^R(\dfrac{2R}{R^2 -|w|^2})^{p}|dw|  < + \infty, $$
where $d_{\Gamma_0}$ denotes the length of the divergent curve $\Gamma_0$ in $M,$ contradicting the assumption of completeness of $M.$ Claim \ref{Cl1} is proved.\\
To summarize, in step 1 we have constructed, for $A = \{z  :  0 < 1/r \leq |z| < r < \infty  \}$
and $S = \{ z : W_z(g_0, g_1)(z) = 0 \}$, a continuous and nowhere vanishing metric 
$d\tau^2$ on $A \setminus S$ which is flat, independent of the choice of coordinate $z$,
and complete with respect to the points of $S$ and with respect to the (outside) boundary
$\{ z : |z| = r\}$.

{\bf Step 2:} We will "symmetrize" the metric constructed in step 2 so that it will become
a complete and flat metric on $Int(A) \setminus (S \cup \tilde{S})$ (with $\tilde{S}$ another discrete subset).

We introduce a new coordinate $\xi (z) :=1/z$ .
By Proposition \ref{W}~a) we have $S= \{ z : W_z(g_0, g_1)(z) = 0 \} = 
\{ z : W_{\xi}(g_0, g_1)(z) = 0 \} $ (where the zeros are taken with the same multiplicities) and since $d\tau^2$ is independent of the coordinate $z$, the change of coordinate $\xi (z) = 1/z$ yields an isometry of $A \setminus S$
onto the set $\tilde{A} \setminus \tilde{S}$, where $\tilde{A}:=\{z : 1/r < |z| \leq r \}$ 
and $\tilde{S}:= \{ z : W_z(g_0, g_1)(1/z) = 0 \}$.
In particular we have (if still $\tilde{h}$ is defined with respect to the coordinate $\xi$) :
$$d\tau^2 = |\tilde{h}(1/z)|^{\frac{2}{1-p}}\bigg(\dfrac{\Pi_{j=1}^q|G_j(1/z)|^{1-\delta}}{e^{\sum_{j=1}^qu_j(1/z)}|W_{(1/z)}(g_0,g_1)(1/z)|}\bigg)^{\frac{2p}{1-p}}|d(1/z)|^2 $$
 $$= \bigg( |h(1/z)|^{\frac{2}{1-p}}\bigg(\dfrac{\Pi_{j=1}^q|G_j(1/z)|^{1-\delta}}{e^{\sum_{j=1}^qu_j(1/z)}|W_z(g_0,g_1)(1/z)|}\bigg)^{\frac{2p}{1-p}} |\frac{dz}{d(1/z)}|^2 \bigg) |d(1/z)|^2$$
 $$= |h(1/z)|^{\frac{2}{1-p}}\bigg(\dfrac{\Pi_{j=1}^q|G_j(1/z)|^{1-\delta}}{e^{\sum_{j=1}^qu_j(1/z)}|W_z(g_0,g_1)(1/z)|}\bigg)^{\frac{2p}{1-p}} |dz|^2$$

\indent We now define 
\begin{align*}
d\tilde{\tau}^2&= \bigg(|h(z)h(1/z)| \cdot \dfrac{\Pi_{j=1}^q|G_j(z)G_j(1/z)|^{(1-\delta)p}}{e^{p\sum_{j=1}^qu_j(z)+p\sum_{j=1}^qu_j(1/z)}|W_z(g_0,g_1)(z) \cdot W_z(g_0,g_1)(1/z)|^p}\bigg)^{\frac{2}{1-p}}|dz|^2\\
&=\lambda^2(z)|dz|^2, 
\end{align*}
on $\tilde{A}_1 := \{z : 1/r < |z| < r \} \setminus \{ z : W_z(g_0, g_1)(z)\cdot W_z(g_0, g_1)(1/z) = 0 \} $. Then $d\tilde{\tau}^2$ is complete  on $\tilde{A}_1$ : In fact by what we showed
above we have:  Towards any point of the boundary 
$\partial \tilde{A}_1 := \{z : 1/r = |z| \} \cup \{z :  |z| = r \} \cup \{ z : W_z(g_0, g_1)(z)\cdot W_z(g_0, g_1)(1/z) = 0 \} $ of $\tilde{A}_1$, one of the factors of $\lambda^2(z)$ is bounded
from below away from zero, and 
the
other factor is the one of a complete metric with respect of this part of the boundary.
Moreover by  the corresponding properties of the two factors of $\lambda^2(z)$ it is trivial that $d\tilde{\tau}^2$ is a continuous nowhere vanishing and flat metric  on $\tilde{A}_1$.

{\bf Step 3:} We produce a contradiction by using Lemma \ref{L5} to the open Riemann surface $(\tilde{A}_1, d\tilde{\tau}^2)$ : \\
In fact, we apply Lemma \ref{L5} to any point $p \in \tilde{A}_1$. Since $d\tilde{\tau}^2$ is 
complete, there cannot exist a divergent curve from $p$ to the boundary $\partial \tilde{A}_1$
with finite length with respect to $d\tilde{\tau}^2$. Since $\Phi : \Delta_{R_0} \rightarrow \tilde{A}_1$ is a local
isometry, we necessarily have $R_0 = \infty$. So $\Phi : {\mathbb C} \rightarrow \tilde{A}_1 \subset \{z : |z| <r\}$ is a non constant holomorphic map, which contradicts to
Liouville's theorem. So our assumption (\ref{ass2}) was wrong.
This proves the Theorem~\ref{T1}. 
\end{proof}
\section{The proof of  Theorem \ref{T2} }
\begin{proof}
 For convenience of the reader, we first recall some notations on the Gauss map of minimal surfaces in $\mathbb R^4$.
Let $x=(x_1, x_2, x_3, x_4) : M \rightarrow \mathbb R^4$ be a non-flat complete minimal surface in $\mathbb R^4.$ As  is well-known, the set of all oriented 2-planes in $\mathbb R^4$ is canonically identified with the quadric  
$$Q_2(\mathbb C):= \{(w_1:...: w_4) | w^2_1 +  ... + w^2_4 = 0\}$$
in $\mathbb P^3(\mathbb C).$ By definition, the Gauss map $g : M \rightarrow Q_2(\mathbb C)$ is the map which maps each point $p$ of $M$ to the point of $Q_2(\mathbb C)$ corresponding to the oriented tangent
plane of $M$ at $p.$ The quadric $Q_2(\mathbb C)$ is biholomorphic to $\mathbb P^1(\mathbb C)\times\mathbb P^1(\mathbb C) .$ By suitable identifications we may regard $g$ as a pair of meromorphic functions $g=(g^1, g^2)$ on $M.$  Let $z$ be a local holomorphic coordinate. 
Set $\phi_i := \partial x_i/dz$ for $i=1,...,4.$ Then, $g^1$ and $g^2$ are given by 
$$g^1 = \dfrac{\phi_3 + \sqrt{-1}\phi_4}{\phi_1 - \sqrt{-1}\phi_2},\  g^2 = \dfrac{-\phi_3 + \sqrt{-1}\phi_4}{\phi_1 - \sqrt{-1}\phi_2}$$
and the metric on $M$ induced from $\mathbb R^4$ is given by
$$ds^2 =|\phi|^2(1 +|g^1|^2)(1+|g^2|^2)|dz|^2 ,$$
where $\phi:= \phi_1 - \sqrt{-1}\phi_2.$
We remark that although the $\phi_i$, $(i=1,2,3,4)$ and $\phi$ depend on $z$, $g=(g^1,g^2)$ and $ds^2$ do not.
Next we take reduced representations $g ^l= (g^l_0 : g^l_1)$ on $M$ and set $||g^l|| = (|g^l_0|^2 +|g^l_1|^2)^{1/2}$ for $l=1,2.$ Then we can rewrite 
\begin{equation}ds^2 = |h|^2||g^1||^2||g^2||^2|dz|^2 \,,
\end{equation}
 where $h:= \phi/(g^1_0g^2_0)$.
 In particular, $h$ is a holomorphic map without zeros. We remark that  $h$ depends on $z$, however, the reduced representations $g^l=(g^l_0:g^l_1)$ are globally defined on $M$ and independent of $z$.
Finally we observe that by the assumption that $M$ is not flat,  $g$ is not constant.

Now the proof of Theorem \ref{T2} will be given in three steps :

{\bf Step 1:} This step is completely analogue to step 1 in the proof of Theorem \ref{T1}.

We may assume that  $ \sum_{j=1}^{q_1}\delta_{g^1}^H(a^{1j})> 2$, $\sum_{j=1}^{q_2}\delta_{g^2}^H(a^{2j})> 2,$ and
\begin{equation}\label{ass4}
\dfrac{1}{\delta_{g^1}^H(a^{1j}) - 2} + \dfrac{1}{\delta_{g^2}^H(a^{2j}) - 2} < 1\,,
\end{equation}
since otherwise case (i) of Theorem \ref{T2} is already proved.

By definition, there exist constants $\eta^l_j \geq 0 (1\leq j \leq q, l=1,2)$ and continous functions $u^l_j (1\leq j \leq q, l=1,2)$ on $M$ such that each $u^l_j$ is harmonic on $M\setminus {(g^l)}^{-1}(a^{lj}), (l=1,2),$ and $u^l_j$ satisfies conditions (D1) and (D2) and $ \gamma_l=q_l-2-\sum_{j=1}^{q_l}\eta^l_j > 0 (l=1,2),$ and
\begin{equation*}
\dfrac{1}{\gamma_1 } + \dfrac{1}{\gamma_2 } < 1\,.
\end{equation*}

Choose  $\delta_0 (> 0)$ such that $\gamma_l - q_l\delta_0 > 0$ for all $l=1,2,$ and
$$\dfrac{1}{\gamma_1 - q_1\delta_0} + \dfrac{1}{\gamma_2 - q_2\delta_0} = 1.$$
If we choose a positive constant $\delta (< \delta_0)$ sufficiently near to $\delta_0$ and set $$p_l := 1/(\gamma_l - q_l\delta),(l=1, 2),$$ we have
\begin{equation} \label{4.2}
0 < p_1 + p_2 < 1, \ \dfrac{\delta p_l}{1-p_1-p_2} > 1 \:( l=1,2)  \ .
\end{equation}
 Consider the  subset
$$ A_2 = A \setminus \{ z : \Pi_{l=1,2}(\Pi_{j=1}^{q_l}|G^l_j(z)|\cdot W_z(g^l_0, g^l_1)(z)) = 0 \}$$
of $A$.  We define a new metric
$$d\tau^2 =  \bigg(|h|\dfrac{\Pi_{j=1}^{q_1}|G^1_j|^{(1-\delta)p_1}\Pi_{j=1}^{q_2}|G^2_j|^{(1-\delta)p_2}}{e^{p_1\sum_{j=1}^{q_1}u^1_j}|W(g^1_0,g^1_1)|^{p_1}e^{p_2\sum_{j=1}^{q_2}u^2_j}|W(g^2_0,g^2_1)|^{p_2}}\bigg)^{\frac{2}{1-p_1 -p_2}}|dz|^2 \ $$
on $A_2$ 
(where again $G^l_j := a^{lj}_0g^l_1 - a^{lj}_1g^l_0\: ( l= 1,2)$ and $h$ is defined with respect
to the coordinate $z$ on $A_2 \subset A$ and $W(g^l_0,g^l_1) = W_z(g^l_0,g^l_1)$).

It is easy to see that by the same arguments as in step 1 of the proof of Theorem \ref{T1}
(applied for each $l=1,2$), we get that $d\tau$ is a continuous nowhere vanishing and flat metric on $A_2$, which is
moreover independant of the choice of the coordinate $z$.

The key point is to prove the following claim :
\begin{claim} \label{Cl2}
 $d\tau^2$ is complete on the set $\{ z : |z| = r\}\cup \{z : \Pi_{l=1,2}W(g^l_0, g^l_1)(z) = 0 \},$ i.e., the set $\{ z : |z| = r\}\cup \{z : \Pi_{l=1,2}W(g^l_0, g^l_1)(z)=0 \}$ is at infinite distance from any interior point in $A_2$. 
\end{claim}
It is easy to see that by the same method as in the proof of Claim~\ref{Cl1} in the proof of Theorem \ref{T1}, we may show that $d\tau$ is complete on $\{z : \Pi_{l=1,2}W(g^l_0, g^l_1)(z) = 0 \}.$

Now assume $d\tau$ is not complete on $\{z : |z| = r \}.$ Then there exists $\rho: [0, 1) \rightarrow A_2,$ where $\rho (1) \in \{z : |z| = r \}$, so that $|\rho| < \infty.$ Furthermore, we may also assume that $dist(\rho (0), \{z : |z| = 1/r\}) > 2|\rho|.$ Consider a small disk $\Delta$ with center at $\rho (0).$ Since $d\tau$ is flat, $\Delta$ is isometric to an ordinary disk in the plane. Let $\Phi: \{|w| < \eta \}\rightarrow \Delta$ be this isometry. Extend $\Phi$, as a local isometry into $A_2,$ to the largest disk $\{|w| < R\} = \Delta_R$ possible. Then $R \leq |\rho|.$ The reason that $\Phi$ cannot be extended to a larger disk is that the image goes to the outside boundary $\{z : |z| = r\}$ of $A_2.$ More precisely, there exists a point $w_0$ with $|w_0| =R$ so that $\Phi(\overline{0,w_0}) = \Gamma_0$ is a divergent curve on $A.$\\
The map $\Phi(w)$ is locally biholomorphic, and the metric on $\Delta_R$ induced from $ds^2$ through $\Phi$ is given by
\begin{equation} \label{4.3} \Phi^*ds^2 =  |h \circ \Phi|^2||g^1 \circ \Phi||^2||g^2 \circ \Phi||^2|\frac{dz}{dw}|^2|dw|^2 \  .
\end{equation}
On the other hand, $\Phi$ is isometric, so we have
$$ |dw| = |d\tau|= \bigg(|h|\dfrac{\Pi_{j=1}^{q_1}|G^1_j|^{(1-\delta)p_1}\Pi_{j=1}^{q_2}|G^2_j|^{(1-\delta)p_2}}{e^{p_1\sum_{j=1}^{q_1}u^1_j}|W(g^1_0,g^1_1)|^{p_1}e^{p_2\sum_{j=1}^{q_2}u^2_j}|W(g^2_0,g^2_1)|^{p_2}}\bigg)^{\frac{1}{1-p_1 -p_2}}|dz|$$
$$\Rightarrow |\dfrac{dw}{dz}|^{1-p_1-p_2} = |h|\dfrac{\Pi_{j=1}^{q_1}|G^1_j|^{(1-\delta)p_1}\Pi_{j=1}^{q_2}|G^2_j|^{(1-\delta)p_2}}{e^{p_1\sum_{j=1}^{q_1}u^1_j}|W(g^1_0,g^1_1)|^{p_1}e^{p_2\sum_{j=1}^{q_2}u^2_j}|W(g^2_0,g^2_1)|^{p_2}}.$$
For each $l=1,2,$ we set $f^l:= g^l(\Phi), f^l_0 := g^l_0(\Phi), f^l_1 := g^l_1(\Phi),  u^l_j:=u^l_j(\Phi)$ and $F^l_j:= G^l_j(\Phi).$  Since
$$W_w(f^l_0, f^l_1) = (W_z(g^l_0, g^l_1) \circ \Phi)\frac{dz}{dw}, (l=1,2), $$
we obtain
\begin{equation} \label{4.4} 
 |\dfrac{dz}{dw}| = \dfrac{\Pi_{l=1,2}(e^{p_l\sum_{j=1}^{q_l}u^l_j}|W(f^l_0,f^l_1)|^{p_l})}{|h(\Phi)|\Pi_{l=1,2}\Pi_{j=1}^{q_l}|F^l_j|^{(1-\delta)p_l}}\ .
 \end{equation}
By (\ref{4.3}) and (\ref{4.4}), we get
\begin{align*}
\Phi^*ds^2& = \bigg(\Pi_{l=1,2} \dfrac{||f^l||(e^{p_l\sum_{j=1}^{q_l}u^l_j}|W(f^l_0,f^l_1)|)^{p_l}}{ \Pi_{j=1}^{q_l}|F^l_j|^{(1-\delta)p_l}}\bigg)^2|dw|^2\\
&= \Pi_{l=1,2}\bigg( \dfrac{||f^l||^{\gamma_l-q_l\delta}e^{\sum_{j=1}^{q_l}u^l_j}|W(f^l_0,f^l_1)|}{ \Pi_{j=1}^q|F^l_j|^{1-\delta}}\bigg)^{2p_l}|dw|^2.
\end{align*}
Using the Lemma \ref{l1}, we obtain
$$\Phi^*ds^2 \leqslant C^{2(p_1+p_2)}.(\dfrac{2R}{R^2 -|w|^2})^{2(p_1+p_2)}|dw|^2.$$
 Since $0 < p_1+p_2 < 1$ by (\ref{4.2}), it then follows that 
$$d_{\Gamma_0} \leqslant \int_{\Gamma_0}ds = \int_{\overline{0,w_0}}\Phi^*ds \leqslant C^{p_1+p_2}. \int_0^R(\dfrac{2R}{R^2 -|w|^2})^{p_1+p_2}|dw|  < + \infty, $$
where $d_{\Gamma_0}$ denotes the length of the divergent curve $\Gamma_0$ in $M,$ contradicting the assumption of completeness of $M.$ Claim \ref{Cl2} is proved.

 {\bf Steps 2 and 3 for the case (i):}
 These steps are  analogue to the corresponding steps in the proof of Theorem \ref{T1}.
Define $d\tilde{\tau}^2 =\lambda^2(z) |dz|^2 \ $ on 
$$\tilde{A}_2 := \{ z : 1/r < |z| < r \}  \setminus $$
$$\setminus \{ z : W_z(g^1_0, g^1_1)(z) \cdot W_z(g^2_0, g^2_1)(z) \cdot W_z(g^1_0, g^1_1)(1/z) \cdot W_z(g^2_0, g^2_1)(1/z)= 0 \} \, , 
$$
 where
\begin{align*}
&\lambda(z)= \bigg(|h(z)|\Pi_{l=1,2}\dfrac{\Pi_{j=1}^{q_l}|G^l_j(z)|^{(1-\delta)p_l}}{e^{p_l\sum_{j=1}^{q_l}u^l_j(z)}|W_z(g^l_0,g^l_1)(z)|^{p_l}}\bigg)^{\frac{1}{1-p_1 -p_2}}\\
& \times \bigg(|h(1/z)|\Pi_{l=1,2}\dfrac{\Pi_{j=1}^{q_l}|G^l_j(1/z)|^{(1-\delta)p_l}}{e^{p_l\sum_{j=1}^{q_l}u^l_j(1/z)}|W_z(g^l_0,g^l_1)(1/z)|^{p_l}}\bigg)^{\frac{1}{1-p_1 -p_2}}.
\end{align*}
By using Claim \ref{Cl2}, the continuous nowhere vanishing and flat metric $d\title{\tau}$
on $A_2$ is also complete.  Using the identical argument of step 3 in the proof of Theorem \ref{T1}
to the open Riemann surface $(\tilde{A}_2, d\tilde{\tau})$ produces a contradiction, so assumption (\ref{ass4}) was wrong.
This implies case (i) of the Theorem \ref{T2}.\\

\indent We finally consider the case (ii) of Theorem \ref{T2} (where $g^2 \equiv constant$ and $g^1 \not\equiv constant$). Suppose that 
$\sum_{j=1}^{q_1}\delta^H_{g^1}(a^{1j}) > 3.$ 

By definition, there exist constants $\eta^1_j \geq 0 (1\leq j \leq q)$ such that $\gamma_1 := q_1-2-\sum_{j=1}^q\eta^1_j > 1$ and continous functions $u^1_j (1\leq j \leq q)$ on $M$ such that each $u^1_j$ is harmonic on $M\setminus (g^1)^{-1}(a^{1j})$ and $u^1_j$ satisfies conditions (D1) and (D2). 

We can choose $\delta$ with 
$$ \dfrac{\gamma_1 - 1}{q_1} > \delta > \dfrac{\gamma_1 - 1}{q_1 +1}, $$
and set $p = 1/ (\gamma_1 -q_1\delta).$ Then 
$$0 < p < 1 , \ \frac{p}{1-p} > \frac{\delta p}{1-p} > 1 . $$
Set $$d\tau^2 = |h|^{\frac{2}{1-p}}\bigg(\dfrac{\Pi_{j=1}^{q_1}|G^1_j|^{1-\delta}}{e^{\sum_{j=1}^{q_1}u^1_j}|W(g^1_0,g^1_1)|}\bigg)^{\frac{2p}{1-p}}|dz|^2 .\ $$
Using this metric, by the analogue arguments as in step 1 to step 3 of the proof of Theorem \ref{T1},  we get the case (ii) of Theorem \ref{T2}.
\end{proof}

{\bf Acknowledgements.} 
The research is partially supported by a NAFOSTED grant of Vietnam.

\vspace{1cm}
{\it Pham Hoang Ha\\
Department of Mathematics, Hanoi National University of Education, 136 XuanThuy str., Hanoi, Vietnam\\
Nguyen Hoang Trang\\
Foreign Language Speacialized School, Hanoi National University, 1 Pham Van Dong str., Hanoi, Vietnam
}

\noindent Emails : ha.ph@hnue.edu.vn, trang1503@gmail.com
\end{document}